\documentclass{article}

\usepackage{arxiv}

\usepackage[utf8]{inputenc} 
\usepackage[T1]{fontenc}    
\usepackage{hyperref}       
\usepackage{url}            
\usepackage{booktabs}       
\usepackage{amsfonts}       
\usepackage{nicefrac}       
\usepackage{microtype}      
\usepackage{lipsum}		
\usepackage{graphicx}
\usepackage{natbib}
\usepackage{doi}
\usepackage{amsmath, amsthm, amssymb}

\newcommand{\norm}[1]{\left\Vert#1\right\Vert}

\newcommand{\R}{\mathbb R}
\newcommand{\N}{\mathbb N}

\newcommand{\intP}{\text{Int}\,P}

\newcommand{\F}{\mathcal{F}}
\newcommand{\D}{\mathcal{D}}
\newcommand{\Fpq}{\mathcal{F}_{p,q}}

\newtheorem{theorem}{Theorem}[section]
\newtheorem{lemma}[theorem]{Lemma}
\newtheorem{proposition}[theorem]{Proposition}
\newtheorem{corollary}[theorem]{Corollary}
\newtheorem{definition}[theorem]{Definition}
\newtheorem{example}[theorem]{Example}

\newtheorem{remark}[theorem]{Remark}

\newcommand\mystyle{\everymath{\displaystyle}}
\mystyle
\title{A common fixed point
theorem for two self-mappings defined on strictly convex probabilistic cone metric space}

\author{ \href{https://orcid.org/0000-0002-3816-5287}{\includegraphics[scale=0.06]{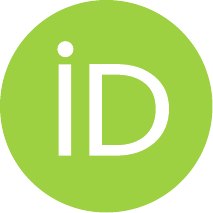}\hspace{1mm}M.H.M.~Rashid}\thanks{Use footnote for providing further
		information about author (webpage, alternative
		address)---\emph{not} for acknowledging funding agencies.} \\
	Department of Mathematics\&Statistics\\Faculty of Science P.O.Box(7)\\
	Mu'tah University University\\
	Mu'tah-Jordan \\
	\texttt{mrash@mutah.edu.jo}
}

\renewcommand{\shorttitle}{\textit{arXiv} Template}

\hypersetup{
pdftitle={A common fixed point
theorem for two self-mappings defined on strictly convex probabilistic cone metric space},
pdfsubject={q-bio.NC, q-bio.QM},
pdfauthor={M.H.M.Rashid},
pdfkeywords={Cone metric space, Normal cones, Non-normal cones, Fixed point theorem},
}

\begin{document}
\maketitle

\begin{abstract}
	This study focuses on defining normal and strictly convex structures within Menger cone PM-space. It also presents a shared fixed point theorem for the existence of two self-mappings constructed on a strictly convex probabilistic cone metric space. The core finding is demonstrated through topological methods to describe spaces with nondeterministic distances. To strengthen our conclusions, we provide several examples. In this research, we introduce and explore normal and strictly convex structures in Menger cone PM-space. A significant contribution of our work is the presentation of a shared fixed point theorem concerning the existence of two self-mappings on a strictly convex probabilistic cone metric space. This theorem is substantiated through topological approaches that effectively describe spaces characterized by nondeterministic distances. To further validate our conclusions, we supplement our theoretical findings with a series of illustrative examples. Our study delves into the intricacies of normal and strictly convex structures within Menger cone PM-space. We present a shared fixed point theorem, demonstrating the existence of two self-mappings in a strictly convex probabilistic cone metric space. Employing topological approaches, we elucidate the key finding and describe spaces with nondeterministic distances. To support and enhance the robustness of our conclusions, we include a variety of examples throughout the study.
\end{abstract}


\section{Introduction}
There are many approaches to generalize metric space. Menger \cite{Menger} created Menger space in 1942 by replacing nonnegative real numbers with distribution functions as metric values. When we do not know the precise distance between two locations but do know the probability of various distance values, we are said to be in a probabilistic metric space.

 A probabilistic generalization of metric spaces appears to be of importance in the study of physical quantities and physiological thresholds. Additionally, it is essential to probabilistic functional analysis. This idea was studied by Schweizer and Sklar  \cite{SS}, while Sehgal and Bharucha-Reid \cite{SBR} were in charge of the significant advancement of Menger space theory. The development of fixed point theory in PM-spaces is credited to Schweizer and Sklar \cite{SS} .\\
\indent  By replacing the set of real numbers with an ordered Banach space, Huang and Zhang in \cite{HZ} extended the idea of metric spaces to create cone metric spaces.

Additionally, they invented the idea of completion in cone metric spaces and provided an explanation of sequence convergence. They established a few fixed point theorems for contractive mappings on whole cone metric space under the presumption of a cone's normality. Since then, a number of researchers have developed Huang and Zhang's discoveries and investigated fixed point theorems for normal and non-normal cones. Numerous research using fixed points have used the Banach contraction idea. This idea has been developed by numerous authors in the form of contraction maps \cite{MSSD, SUH}.

 In contrast, Takahashi\cite{Takahashi} defined the idea of a convex structure in a metric space and referred to such a space as a convex metric space. In addition, he looked at a number of the space's characteristics and made sure that there was a fixed point for non-expansive maps in the context of convex metric space. Many fixed point and common fixed point theorems in the setting of convex metric space have been produced over the past forty years; for examples, see
see \cite{GSW, HL, JNB, Rashid1}.

The aim of this study is to define normal and strictly convex structures within Menger cone PM-space and to present a shared fixed point theorem for the existence of two self-mappings on a strictly convex probabilistic cone metric space. The research employs topological approaches to elucidate the key finding, which describes spaces with nondeterministic distances. To reinforce our conclusions, we provide several illustrative examples. This comprehensive approach not only broadens the understanding of these mathematical structures but also demonstrates their applicability through concrete instances, thereby enhancing the robustness and validity of our theoretical results.

The organization of this paper is as follows: In the next section, we define cone metric spaces and verify several foundational facts. Section three is dedicated to establishing topological properties pertinent to our study. Section four focuses on defining strictly convex and normal structures in probabilistic cone metric spaces. In section five, we demonstrate that two self-mappings on a strictly convex probabilistic cone metric space share a common fixed point theorem. Topological approaches for characterizing spaces with nondeterministic distances are utilized to prove the primary conclusion. Each section builds upon the previous one, creating a comprehensive framework for understanding and validating the main findings of the study. This structured approach ensures clarity and logical progression throughout the paper.

\section{Preliminaries }
We'll define cone metric spaces and verify several facts in this section.
Let $E$ always be a real Banach space, and $P$ always be a subset of $E$. $P$ is a cone if and only if:
\begin{enumerate}
  \item [(P1)] $P$ is closed, nonempty, and $P\neq \{0\}$;
  \item [(P2)] $a, b \in\R$, $a,b \geq 0$, $x,y \in P$ implies  $ax+by\in P;$
  \item [(P3)] $x\in P$ and $-x\in P$ implies $x=0$.
\end{enumerate}
Given a cone $P\subset E$,  we define a partial ordering $\preceq$ with regard to $P$ by $x\preceq y$
 if and only if $y-x \in P$. We will write $x\prec y$  to show that $x preceq y$   but $x\neq y$, while $x\ll y$ will
stand for $y-x\in \text{Int}\,P$ , $\intP$ denotes the interior of $P$.\\
If there is an integer $K>0$ such that for any $x,y\in E$, $$0\leq x\leq y \Longrightarrow \norm{x}\leq K\norm{y},$$ the cone $P$ is said to as normal.
The normal constant of $P$  (see \cite{HZ}) is the smallest positive number that satisfies the criteria above.
Rezapour and Hamlbarani \cite{RH} demonstrated that there are cones with normal constants $K=1$, $M>K$ for any $K>1$, and that there are no cones with normal constants $K<1$.
Every increasing sequence that is bounded from above must converge in order for the cone $P$ to be said to as regular. In other words, if $\{x_n\}$ is a sequence such that $$x_1\leq x_2\leq\cdots\leq x_n\cdots\leq y$$ for some $y \in E$, then $x \in E$ is such that $\norm{x_n-x}\rightarrow 0$ ($n\rightarrow \infty$). In other words, the cone $P$ is regular if and only if each decreasing sequence that is bounded from below is convergent. Regular cones are recognized as being normal cones. \\
\indent In the following, $E$ is always assumed to be a Banach space, $P$ to be a cone in $E$
with $\intP\neq \emptyset$, and $\leq$ to be partial ordering with regard to $P$.

\begin{proposition}(\cite{RH, TA}).\label{P1}
  Let $P$ be a cone of $E.$ Then
  \begin{enumerate}
    \item [(a)] $P+\intP\subset \intP$;
    \item [(b)] For every $\alpha\in \R^+$, we have $\alpha \intP\subseteq\intP$;
    \item [(c)] For each $\theta\leq c_1$ and $\theta\leq c_2$, there is an element $\theta\leq c$  such that $c\leq c_1$, $c\leq c_2$.
  \end{enumerate}
\end{proposition}

\begin{definition}(\cite{HZ}).
  A cone metric space is an ordered $(X,d)$, where $X$ is any set and $d:X\times X\rightarrow E$ is
a mapping satisfying:
\begin{enumerate}
  \item [(CM1)] $d(x,y)\geq \theta$ for all $x,y\in X$,
  \item [(CM2)] $d(x,y)=\theta$ if and only if $x=y$,
  \item [(CM3)] $d(x,y)=d(y,x)$ for all $x,y\in X$,
  \item [(CM4)] $d(x,y)\leq d(x,z)+d(z,y)$ for all $x,y,z\in X$.
\end{enumerate}
\end{definition}

In \cite{TA}, for $c\in E$  with $c\gg \theta$ and $x\in X$, define
$B(x,c)= \{y \in X: d(x,y)\ll c\}$
and
$\beta=\{B(x,c):x\in X,c\in E\,\,\text{with}\,\, c\gg \theta\},$
 then show that
$$\tau_c=\{U\subset X:\forall x\in U, \exists B(x,c)\in\beta,x\in B(x,c)\subset U\}$$
is a topology on $X.$

\begin{definition} (\cite{HZ}).
  Let $(X,d)$ be a cone metric space, $x\in X$  and $\{x_n\}$ be a sequence in $X.$ Then
  \begin{enumerate}
    \item [(i)] The sequence $\{x_n\}$ is said to converge to $x$ if for any $c\in E$  with $c\gg \theta$ there exists a natural number $n_0$ such that
$d(x_n,x)\ll c$ for all $n\geq n_0$. We denote this by $\displaystyle{\lim_{n\rightarrow \infty} x_n=x}$ or $x_n\rightarrow x$ as $n\rightarrow \infty.$
    \item [(ii)] The sequence $\{x_n\}$ is said to be a Cauchy sequence if for any $c\in E$  with $c\gg \theta$ there exists a natural number $n_0$ such that $d(x_n,x_m)\ll c$ for all $n,m\geq n_0$.
    \item [(iii)] $(X,d)$ is said to be a complete cone metric space if every Cauchy sequence is convergent.
  \end{enumerate}
\end{definition}

\begin{definition}(\cite{SS}).
 A binary operation $\ast: [0,1]\times [0,1]\rightarrow [0,1]$ is said to be a
continuous $t$-norm if $([0,1],\ast)$ is a topological monoid with unit 1 such that
$a\ast b\leq c\ast d$  whenever $a\leq c,b\leq d$ for all $a,b,c,d\in [0,1]$.
\end{definition}
\begin{remark}(\cite{GV}).
  For any $r_1>r_2$, we can find a $r_3$ such that $r_1\ast r_3\geq r_2$ and for any $r_4$ we can find a $r_5$ such that
  $r_5\ast r_5\geq r_4,$ where $r_j\in (0,1)$ for $j=1,2,\cdots,5.$
\end{remark}
Some typical examples of $t$-norm are the following:
\begin{align*}
  a\ast b &=ab,\quad \text{(product)} \\
  a\ast b &=\min\{a,b\},\quad \text{(minimum)}
\end{align*}
\section{Probabilistic Cone Metric Spaces}
Let $P$ be a cone and $E$ be a real, ordered Banach space. If a mapping $\F: E\rightarrow int(P)$  is nondecreasing and left continuous with $\inf \F=\theta$ and $\sup \F=1$, it is referred to as a distribution function. The set of all distribution functions will be denoted as $\D$.
\begin{definition}
  If $X$ is an arbitrary set, $\ast$ is a continuous $t$-norm, and $\F$ is a mapping of $X\times X$ into $D$, which associates a distribution function $\F_(p,q)$ with each pair of $(p,q)$ points in $X$, then the triple $(X,\F,\ast)$ is said to be a probabilistic cone metric space The distribution function $\F_(p,q)$ will be denoted by $\F_{pq}$,
where the value of $\F_{pq}$ for the argument $t$ is indicated by the symbol $\F_{pq}(t)$. The axioms below are taken to be true for the function $\F_{pq}$.
  For all $p,q,r\in X$ and $t,s\in int(P)$ (that is, $t\gg \theta,s\gg\theta$)
  \begin{enumerate}
    \item [(PCM1)] $\Fpq(t)>0$;
    \item [(PCM2)] $\Fpq(t)=1$ if and only if $x=y$;
    \item [(PCM3)] $\Fpq=\F_{q,p}$;
    \item [(PCM4)] If $\Fpq(t)=1$ and $\F_{q,r}(s)=1$, then $\F_{p,r}(t+s)=1$.
  \end{enumerate}
\end{definition}
\begin{definition}
   A triple $(X,\F,\ast)$  is said to be a Menger probabilistic cone metric space if
    $(X,\F,\ast)$  is a probabilistic cone metric space and satisfying the following condition:
    \begin{enumerate}
      \item [(PCM5)] $\F_{p,r}(t+s)\geq \Fpq(t)\ast \F_{q,r}(s)$, for all $p,q,r\in X$ and $t,s\in int(P)$.
    \end{enumerate}
\end{definition}
\begin{remark} (i) If we take $E=\R$, $P=\R^+$, then every probabilistic  metric space become a
probabilistic cone metric space.\\
(ii) Every cone metric space may be regarded as a $PCM$-space of a special if we have
only to set $\Fpq(t)=H(t-d(p,q))$ for every pair of points $(p,q)$ in the cone metric space, where
$H$ denote the specific distribution function defined by
$$H(t)=\left\{
         \begin{array}{ll}
           \theta, & \hbox{if $t\leq \theta$;} \\
           1, & \hbox{if $t\gg \theta$.}
         \end{array}
       \right.
$$
\end{remark}
\begin{example}(\cite{Rashid1})
	Let $E=\mathbb{R}^2$. Then $P=\{(x,y)\in\R^2: x,y\geq 0\}$ is a normal cone with normal constant $K=1$.
	Let $X=\mathbb{R}$, $\ast(a,b)=ab$ and $\F:X^2\times int(P)\rightarrow [0,1]$ defined by
	$$\F_{x,y}(t)=\dfrac{1}{e^{\frac{|x-y|}{\norm{t}}}}$$
	for all $x,y\in X$ and $t\gg \theta$. Then $(X,\F,T)$ is a probabilistic cone metric spaces.
\end{example}
\begin{example}(\cite{Rashid1})
	Let $P$ be any cone, $X=\mathbb{N}$, $\ast(a,b)=ab$, $\F:X^2\times int(P)\rightarrow [0,1]$ defined by
	$$\F_{x,y}(t)=\left\{
	\begin{array}{ll}
	x/y, & \hbox{if $x\leq y$;} \\
	y/x, & \hbox{if $y\leq x$.}
	\end{array}
	\right.
	$$
	for all $x,y\in X$ and $t\gg \theta$. Then $(X,\F,T)$ is a probabilistic cone metric spaces.
\end{example}
\begin{definition}
  Let $p$ be a point in a $PCM$-space $(X,\F)$. By an $(\epsilon,\lambda)$-neighborhood of
  $p$, $\epsilon\gg \theta$,$\lambda\in (0,1)$, we mean the set of points $q$ in $X$ for which
  $\Fpq(\epsilon)>1-\lambda$, that is,
  $$N_{p}(\epsilon,\lambda)=\{q\in X: \Fpq(\epsilon)>1-\lambda\}.$$
\end{definition}
\begin{proposition}
  If $\epsilon_1\leq \epsilon_2$ and $\lambda_1\leq \lambda_2$, then $N_p(\epsilon_1,\lambda_1)\subset
  N_p(\epsilon_2,\lambda_2)$.
\end{proposition}
\begin{proof}
  Suppose $q\in N_p(\epsilon_1,\lambda_1)$ so that $\Fpq(\epsilon_1)>1-\lambda_1$. Then
  we have $\Fpq(\epsilon_2)\geq \Fpq(\epsilon_1)>1-\lambda_1\geq 1-\lambda_2$ and so by definition,
  $q\in  N_p(\epsilon_2,\lambda_2)$.
\end{proof}
\begin{proposition}
  Let $(X,\F,\ast)$ be a Menger probabilistic cone metric space. Define
  $$\tau_P=\{A\subset X:x\in A\Longleftrightarrow \exists  \epsilon\gg \theta, \lambda\in (0,1), \,\,\text{such that}\,\, N_p(\epsilon,\lambda)\subset A\}.$$
  Then $\tau_P$ is a topology on $X$.
\end{proposition}
\begin{proof}(i) $\emptyset\in \tau_{P}$ since $\emptyset$ has no
elements. Since for any $p\in X$, any $\lambda\in (0,1)$ and any $\epsilon\gg \theta$, $N_p(\epsilon,\lambda)\subset X$, then $X\in \tau_{P}$.\\
(ii) Let $A_i\in\tau_P$ for each $i\in I$ and $p\in \bigcup_{i\in I}A_i$. Then there exists $i_0\in I$ such that $p\in A_{i_0}$.
So, there exist $\epsilon\gg \theta$ and $\lambda\in (0,1)$ such that $N_p(\epsilon,\lambda)\subset A_{i_0}\subset \bigcup_{i\in I}A_i$
and hence $\bigcup_{i\in I}A_i\in\tau_P$.\\
(iii) Let $A_1,A_2,\cdots,A_n\in\tau_{P}$ and $p\in \bigcap_{i=1}^{n} A_i$. Then $p\in A_i$, for $i=1,\cdots,n$, so there exists $\epsilon_i\gg \theta, i=1,\cdots,n$, and $\lambda_i\in (0,1),i=1,\cdots,n$ such that $N_p(\epsilon_i,\lambda_i)\subset A_i$  for $i=1,\cdots,n$. It follows from Proposition \ref{P1}
that, for $\epsilon_i\gg \theta,i=1,\cdots,n$, there exists $\epsilon\gg \theta$ such that $\epsilon\ll \epsilon_i,i=1,\cdots,n$ and take
$\lambda=\min\{\lambda_1,\lambda_2,\cdots,\lambda_n\}$. Hence $N_p(\epsilon,\lambda)\subset \bigcap_{i=1}^{n}N_p(\epsilon_i,\lambda_i)\subset \bigcap_{i=1}^{n}A_i$. Thus $\bigcap_{i=1}^{n}A_i\in \tau_{P}$.
This proves that $\tau_{P}$ is a topology on $X$.
\end{proof}
\begin{proposition}\label{Okasha2}
   Let $(X,\F,\ast)$ be a Menger probabilistic cone metric space. Then $(X,\tau_{P})$ is Hausdorff.
\end{proposition}
\begin{proof}
  Let $p,q\in X$ be two distinct points. Then $0<\F_{pq}(t)<1$. Set $\F_{pq}(t)=\lambda$, for some $\lambda\in (0,1)$. For
  each $\lambda_0$ such that $1>\lambda_0>\lambda$, we can find $\lambda_1\in (0,1)$ such that  $\lambda_1\ast \lambda_1\geq \lambda_0$.\\
  \indent Consider the open balls $N_p\left(\frac{\epsilon}{2},1-\lambda_1\right)$ and $N_q\left(\frac{\epsilon}{2},1-\lambda_1\right)$. Then
  $$N_p\left(\frac{\epsilon}{2},1-\lambda_1\right)\cap N_q\left(\frac{\epsilon}{2},1-\lambda_1\right)=\emptyset.$$
  Otherwise, if $N_p\left(\frac{\epsilon}{2},1-\lambda_1\right)\cap N_q\left(\frac{\epsilon}{2},1-\lambda_1\right)\neq\emptyset$, then there exists a $r\in N_p\left(\frac{\epsilon}{2},1-\lambda_1\right)\cap N_q\left(\frac{\epsilon}{2},1-\lambda_1\right)$. Then $\F_{p,r}(t)>1-(1-\lambda_1)=\lambda_1$ and $\F_{qr}>1-(1-\lambda_1)=\lambda_1$. Thus
  \begin{align*}
    \lambda &=\Fpq(t) \\
     & \geq \F_{p,r}\left(\frac{t}{2}\right)\ast \F_{r,q}\left(\frac{t}{2}\right)\\
     &\geq \lambda_1\ast \lambda_1\geq \lambda_0\\
     &>\lambda
  \end{align*}
  which is a contradiction. Therefore $(X,\F,\ast)$ is a Hausdorff.
\end{proof}
\begin{definition}
  Let $(X,\F,\ast)$ be a Menger probabilistic cone metric space, $x\in X$ and $\{x_n\}$ be a sequence in $X$. Then
  \begin{enumerate}
    \item [(a)] $\{x_n\}$ is said to converge to $x$ if for any $\epsilon\gg \theta$ and any $\lambda\in (0,1)$ there exists a natural number $n_0$ such that
    $\F_{x_n,x}(\epsilon)>1-\lambda$ for all $n\geq n_0$.
     We denote this by $\displaystyle{\lim_{n\rightarrow \infty}x_n=x}$ or $x_n\rightarrow x$ as $n\rightarrow \infty$.
    \item [(b)] $\{x_n\}$ is said to be a Cauchy sequence if for any $\epsilon\gg \theta$ and any $\lambda\in (0,1)$ there exists a natural number $n_0$ such that
    $\F_{x_n,x_m}(\epsilon)>1-\lambda$ for all $n,m\geq n_0$.
    \item [(c)] $(X,\F,\ast)$ is said to be a complete Menger probabilistic cone metric space if every Cauchy sequence is convergent.
  \end{enumerate}
\end{definition}
\begin{proposition}\label{Okasha3}
 Let $(X,\F,\ast)$ be a Menger probabilistic cone metric space. Then $(X,\tau_{P})$ is first countable.
\end{proposition}
\begin{proof}
  Let $\epsilon\gg \theta$, $x\in X$. We will show that $B_x=\{N_x\left(\frac{1}{n},\frac{\epsilon}{n}\right):n\in\N\}$ is a local
  basis for $x\in X$. Let $V\in \tau_{P}$ and $x\in V$. Since $V$ is open, then there exists $\lambda\in (0,1)$ and $\epsilon\gg \theta$ such that
  $N_x(\epsilon,\lambda)\subset V$. Choose $n\in\N$ such that $\frac{1}{n}<\lambda$ and $\frac{\epsilon}{n}\ll \epsilon$. Now we just need to show
  $N_x\left(\frac{\epsilon}{n},\frac{1}{n}\right)\subset N_x(\epsilon,\lambda)$. Let $z\in N_x\left(\frac{\epsilon}{n},\frac{1}{n}\right)$.
  Then $\F_{x,z}\left(\frac{\epsilon}{n}\right)>1-\frac{1}{n}>1-\lambda$.
  Since $\frac{\epsilon}{n}\leq \epsilon,$ we have $1-\lambda<\F_{x,z}\left(\frac{\epsilon}{n}\right)\leq \F_{x,z}(\epsilon)$. Hence $z\in N_x(\epsilon,\lambda)$ which implies that
  $N_x\left(\frac{\epsilon}{n},\frac{1}{n}\right)\subset N_x(\epsilon,\lambda)\subset V$. Consequently, $B_x$ is a countable local basis for $x$.
  Therefore, $(X,\tau_{P})$ is first countable.
\end{proof}
\begin{definition}
 Let $(X,\F,\ast)$ be a Menger probabilistic cone metric space. A subset $A$ of $X$ is said to be $FC$-bounded if there
  exists $\epsilon\gg \theta$ and $\lambda\in (0,1)$ such that $\F_{x,y}(\epsilon)>1-\lambda$ for all $x,y\in A$.
\end{definition}
\begin{proposition}\label{Okasha4}
   Let $(X,\F,\ast)$ be a Menger probabilistic cone metric space. If $S$ is compact subset of $X$, then it is closed and $FC$-bounded.
\end{proposition}
\begin{proof}
  Let $\epsilon\gg \theta$ and $\lambda\in (0,1)$. Let $\mathcal{U}=\{N_x(\epsilon,\lambda):x\in S\}$
 be an open cover of $S$. Since $S$ is compact, we have a subcover of $\mathcal{U}$, i.e.,
 there exist $x_1,x_2,\cdots,x_n\in A$ such that $A\subset \bigcup_{i=1}^{n}N_{x_i}(\epsilon,\lambda)$.
 For any $x,y\in S$ there exist $1\leq i,j\leq n$ such that $x\in N_{x_i}(\epsilon,\lambda)$ and $y\in N_{x_j}(\epsilon,\lambda)$.
 Therefore $\F_{x,x_i}(\epsilon)>1-\lambda$ and $\F_{y,x_j}(\epsilon)>1-\lambda$. Now let $a=\min\{\F_{x_i,x_j}(\epsilon):1\leq i,j\leq n\}$.
 Then $a>0$. Now
 \begin{align*}
   \F_{x,y}(3\epsilon) & \geq \F_{x,x_i}(\epsilon)\ast \F_{x_i,x_j}(\epsilon)\ast \F_{y,x_j}(\epsilon)\\
    &\geq (1-\lambda)\ast a\ast (1-\lambda)
 \end{align*}
  Let $\epsilon'=3\epsilon$ and choose $0<\mu<1$ such that $(1-\lambda)\ast a\ast (1-\lambda)>1-\mu, 0<\mu<1$,
 we have $\F_{x,y}(\epsilon')>1-\mu$ for all $x,y\in S$.  Hence $S$ is $FC$-bounded. On the other hand, since a Menger probabilistic cone metric space is Hausdorff and every compact subset of a Hausdorff space is closed, $S$ is closed.
\end{proof}
\begin{proposition}\label{Okasha5}
  Let $(X,\F,\ast)$ be a Menger probabilistic cone metric space and $\tau_P$ be the topology induced by probabilistic cone metric space.
Then for any nonempty subset $S\subseteq X$ we have
\begin{enumerate}
  \item [(i)] $S$ is closed if and only if for any sequence $\{x_n\}$ in $X$ which converges to $x$, we have $x\in S$;
  \item [(ii)] if we define $\overline{S}$ to be the intersection of all closed subset of $X$ which contain $S$, then for any $x\in \overline{S}$ and for any $0<\lambda <1$ and $\epsilon\gg \theta$, we have $N_x(\epsilon,\lambda)\cap S\neq \emptyset$.
\end{enumerate}
\end{proposition}
\begin{proof} (i) Assume that $S$ is closed and let $\{x_n\}$ be a sequence in $S$ such that $\displaystyle{\lim_{n\rightarrow \infty}x_n=x}$. Let us prove that
$x\in S$. Assume not, i.e. $x\notin S$. Since $S$ is closed, then there exists $0<\lambda <1$ and $\epsilon\gg \theta$ such that $N_x(\epsilon,\lambda)\cap X=\emptyset$.
Since $\{x_n\}$ converges to $x$, then there exists $N\geq 1$ such that for any $n\geq N$ we have $x_n\in N_x(\epsilon,\lambda)$. Hence
$x_n\in N_x(\epsilon,\lambda)\cap S$, which leads to a contradiction. Conversely assume that for any sequence $\{x_n\}$ in $S$ which
converges to $x$, we have $x\in S$. Let us prove that $S$ is closed. Let $x\notin S$. We need to prove that there
exists $0<\lambda<1$  and $\epsilon\gg \theta$ such that $N_x(\epsilon,\lambda)\cap S=\emptyset$. Assume not, i.e. for any $0<\lambda <1$ and $\epsilon\gg \theta$, we
have $N_x(\epsilon,\lambda)\cap S\neq \emptyset$. So for any $n\geq 1$, choose $x_n\in N_x(\epsilon, \frac{1}{n})$. Clearly we have $\{x_n\}$ converges to $x$. Our
assumption on $S$ implies $x\in S$,  a contradiction.\\
(ii) Clearly $\overline{S}$ is the smallest closed subset which contains $S$. Set
$$S^*=\{x\in X:\text{for any $\epsilon\gg \theta$, there exists $a\in S$ such that $\F_{x,a}(\epsilon)>1-\lambda$}\}.$$
We have $S\subset S^*$. Next we prove that $S^*$ is closed. For this we use property (i). Let $\{x_n\}$  be a sequence in $S^*$
such that $\{x_n\}$ converges to $x$. Let $0<\lambda<1$ and $\epsilon\gg \theta$. Since $\{x_n\}$ converges to $x$, there exists $N\geq 1$ such that
for any $n\geq N$ we have
$$\F_{x,x_N} \left(\frac{t}{2}\right)>1-\lambda.$$
 Let $\lambda_0=\F_{x,x_N} \left(\frac{t}{2}\right)>1-\lambda$. Since $\lambda_0> 1-\lambda$, we can find an $\mu$, $0<\mu<1$, such that
$\lambda_0>1-\mu>1-\lambda_0$. Now for a given $\lambda_0$ and $\mu$ such that $\lambda_0>1-\mu$ we can find $\lambda_1$, $0 <\lambda_1 < 1$, such that
$$\lambda_0\ast (1-\lambda_1)\geq 1-\mu.$$
Now since $x_n\in S^*$, there exists $a\in X$ such that
$$ \F_{x_n,a}\left(\frac{\epsilon}{2}\right)>1-\lambda_1$$
Hence
$$\F_{x,a}(\epsilon)\geq \F_{x,x_N} \left(\frac{t}{2}\right)\ast \F_{x_n,a} \left(\frac{\epsilon}{2}\right)>\lambda_0\ast (1-\lambda_1)\geq 1-\mu>1-\lambda,$$
which implies $x\in S^*$. Therefore $S^*$ is closed and contains $S$. The definition of $\overline{S}\subset S^*$, which implies the
conclusion of (ii).
\end{proof}
Note that, every compact subset of a Hausdorff topological space is closed.
\begin{proposition}\label{Okasha6}
  Let $(X,\F,\ast)$ be  a Menger probabilistic cone metric space and $\tau_M$ be the topology induced by probabilistic cone metric type.
Let $S$ be a nonempty subset of $X$. The following properties are equivalent:
\begin{enumerate}
  \item [(i)] $S$ is compact.
  \item [(ii)] For any sequence $\{x_n\}$ in $S$, there exists a subsequence $\{x_{n_k}\}$ of $\{x_n\}$ which
   converges, and if $\{x_{n_k}\}$ converges to $x$ then $x\in S$.
\end{enumerate}
\end{proposition}
\begin{proof}
  (i) Assume that $S$ is a nonempty compact subset of $X$. It is easy to see that any decreasing sequence of
nonempty closed subsets of $S$ have a nonempty intersection. Let $\{x_n\}$ be a sequence in $S$. Set
$C_m=\{x_m:m\geq n\}$. Then we have $\bigcap_{n\geq 1}\overline{C_n}\neq \emptyset$. Then for $0<\lambda< 1$, $t\gg \theta$ and for any $n\geq 1$, there
exists $m_n\geq n$ such that $\F_{x,x_{m_n}}(\epsilon)>1-\lambda$. This clearly implies the existence of a subsequence of $\{x_n\}$ which
converges to $x$. Since $S$ is closed, then we must have $x\in S$.\\
Conversely let $S$ be a nonempty subset of $X$ such that the conclusion of (ii) is true. Let us prove that $S$ is
compact. First note that for any $0<\lambda<1$, $t\gg \theta$, there exists $x_1,x_2,\cdots,x_n\in A$ such that
$$S\subseteq \bigcup_{i=1}^{n}N_{x_i}(\epsilon,\lambda).$$
Assume not, then there exists $0 <\lambda_0< 1$, such that for any finite number of points $x_1,x_2,\cdots,x_n\in X$, we have
$$S\nsubseteq \bigcup_{i=1}^{n}N_{x_i}(\epsilon,\lambda).$$
Fix $x_1\in X$. Since $S\nsubseteq N_{x_1}(\epsilon,\lambda)$, there exists $x_2\in S\setminus N_{x_1}(\epsilon,\lambda)$. By induction we build a sequence $\{x_n\}$ such that
$$x_{n+1}\in S\setminus (N_{x_1}(\epsilon,\lambda)\cup\cdots \cup N_{x_n}(\epsilon,\lambda))$$
for all $n\geq 1$. Clearly we have $\F_{x_n,x_m}(\epsilon) <1-\lambda_0$ for all $n,m\geq 1$, with $n\neq m$. This condition implies
that no subsequence of $\{x_n\}$ will be Cauchy or convergent. This contradicts our assumption on $X$. Next let
$\{O_{\alpha}\}_{\alpha\in J}$ be an open cover of $S$. Let us prove that only finitely many $O_{\alpha}$ cover $S$. Fix $\epsilon\gg \theta$, first note that there exists $0<\lambda_0< 1$ such that for any $x\in S$, there exists $\alpha\in J$ such that $N_x(\lambda_0,\epsilon)\subset O_{\alpha}$. Assume not, then for
any $0 <\lambda< 1$, there exists $x_{\lambda}\in X$ such that for any $\alpha\in J$, we have $N_x(\epsilon,\lambda)\nsubseteq O_{\alpha}$. In particular, for any $n\geq 1$,
there exists $x_n\in X$ such that for any $\alpha\in J$, we have $N_x(\epsilon,1/n)\nsubseteq O_{\alpha}$. By our assumption on $S$, there exists
a subsequence $\{x_{n_k}\}$ of $\{x_n\}$ which converges to some point $x\in X$. Since the family $\{O_{\alpha}\}_{\alpha\in J}$ covers $X$, there
exists $\alpha_0\in J$ such that $x\in O_{\alpha_0}$. Since $O_{\alpha_0}$ is open, there exists $0<\lambda_0< 1$, and $\epsilon_0\gg \theta$ such that $N_x(\epsilon,\lambda)\subset O_{\alpha_0}$ .
Fix $\epsilon\gg \theta$ and let $\epsilon_1=\epsilon$, for any $n_1\geq 1$ and $a\in N_{n_{n_1}}(\epsilon,\frac{1}{n_1})=N_{n_{n_1}}(\epsilon_1,\frac{1}{n_1})$ we have
$$\F_{x,a}(\epsilon_0)\geq \F_{x,x_{n_k}}(\epsilon_0-\epsilon_1)\ast \F_{x_{n_k},a}(\epsilon_1)
 >\F_{x,x_{n_k}}(\epsilon_0-\epsilon_1)>1-\frac{1}{n_k}$$
for $n_k$ large enough, we will get $\F_{x,a}(\epsilon)>1-\lambda_0$ for any $a\in N_{x_{n_k}}(\epsilon,\frac{1}{n_k})$. In the other words, we have
$B_{x_{n_k}}(\frac{1}{n_k},t)\subset B_x(r_0,t_0$, which implies
$$N_{x_{n_k}}(\epsilon,\frac{1}{n_k})\subset O_{\alpha_0}.$$
This is in clear contradiction with the way the sequence $\{x_n\}$ was constructed. Therefore there exists
$0<\lambda_0<1$ such that for any $x\in S$, there exists $\alpha\in J$ such that $N_x(\epsilon,\lambda_0)\subset O_{\alpha}$. For such $\lambda_0$, there exist
$x_1,x_2,\cdots,x_n\in X$ such that
$$S\subset N_{x_1}(\epsilon,\lambda_0)\cup\cdots\cup N_{x_n}(\epsilon,\lambda_0)$$
But for any $i=1,2,\cdots,n$ there exists $\alpha\in J$ such that $N_{x_i}(\epsilon,\lambda_0)\subset O_{\alpha}$, i.e. $S\subset O_{\alpha_1}\cup\cdots\cup O_{\alpha_n}$. This completes the
proof that $S$ is compact.
\end{proof}
\begin{definition}
   Let $(X,\F,\ast)$ be  a Menger probabilistic cone metric space. The subset $S$ of $X$ is called sequentially compact
if and only if for any sequence $\{x_n\}$ in $S$, there exists a subsequence $\{x_{n_k}\}$
 of $\{x_n\}$ which converges to $x$, and $x\in X$. Also $S$ is called totally bounded if for any $0<\lambda<1$ and $\epsilon\gg \theta$ there exist
  $x_1,x_2,\cdots,x_n\in X$ such that $S\subset \bigcup_{i=1}^{n}N_{x_i}(\epsilon,\lambda)$.
\end{definition}
The following result is a simple consequence of the last results about Menger probabilistic cone metric space.
\begin{corollary}
  Let $(X,\F,\ast)$ be  a Menger probabilistic cone metric space. Then
  \begin{enumerate}
    \item [(i)] Every compact set $S$ of $X$ is complete.
    \item [(ii)] Every closed subset of a complete probabilistic metric type space is complete.
  \end{enumerate}
\end{corollary}
\begin{definition}
Let $(X,\F,\ast)$ be a probabilistic cone metric space with a continuous $t$-norm $\ast$.
Let $\lambda\in (0,1)$, $\epsilon\gg \theta$ and $x\in X$. The set $N_{x}[\epsilon,\lambda]=\{y\in X:\F_{x,y}(\epsilon)\geq 1-\lambda\}$
is called closed $(\epsilon,\lambda)$-neighborhood of a point $x$ in $X$.
\end{definition}
\begin{definition}
  A subset $K$ of a probabilistic cone metric space is called compact if the following statement holds:
  $$K\subseteq \displaystyle{\bigcup_{\alpha\in \Lambda} U_{\alpha}}\Longrightarrow
  K\subseteq \displaystyle{\bigcup_{i=1}^{n}U_{\alpha_i}}
  \,\,\text{for some $\alpha_1,\cdots,\alpha_n\in \Lambda$ }$$
  for every collection $\{U_{\alpha}:\alpha\in\Lambda\}$ of open sets $U_{\alpha}\subseteq X$.
\end{definition}
\begin{lemma}\label{Lemma2.9}
  Let $(X,\F,\ast)$ be a probabilistic cone metric space with a continuous $t$-norm $\ast$ and $K\subseteq X$.
  Then, $K$ is compact if and only if for every collection of closed sets $\{F_{\alpha}\}$ such that
  $F_{\alpha}\subseteq K$ it holds that
  $$\displaystyle{\bigcap_{\alpha\in \Lambda} F_{\alpha}}\Longrightarrow \displaystyle{\bigcap_{i=1}^{n}F_{\alpha_i}=\emptyset}
   \,\,\text{for some $\alpha_1,\cdots,\alpha_n\in \Lambda$ }.$$
\end{lemma}
\begin{lemma} Let $(X,\F,\ast)$ be a probabilistic cone metric space with a continuous $t$-norm $\ast$ and $K\subseteq X$.
Then $x\in \overline{K}$ if and only if there exists a sequence $\{x_n\}$ in $K$ such that $x_n\rightarrow x$.
\end{lemma}
\begin{definition}
   Let $(X,\F,\ast)$ be a probabilistic cone metric space with a continuous $t$-norm $\ast$ and $A\subseteq X$.
   The probabilistic diameter of $A$ is given by
   $$\delta_A(t)=\displaystyle{\inf_{x,y\in A}\sup_{s<t}\F_{x,y}(s)}.$$
   The diameter of the set $A$ is defined as
   $$\delta_A=\displaystyle{\sup_{t>0}\inf_{x,y\in A}\sup_{s<t}\F_{x,y}(s)}.$$
   If there exists a number $\lambda\in (0,1)$ such that $\delta_A=1-\lambda$, then the set $A$
   is called probabilistic semi-bounded. If $\delta_A=1$, then $A$ is called probabilistic bounded.
\end{definition}
\begin{theorem}\label{Theorem2.13}
  Every compact subset $A$ of a probabilistic cone metric space $(X,\F,\ast)$  with continuous
  $t$-norm $\ast$ is probabilistic semi-bounded.
\end{theorem}
\begin{proof}
  Let $A$ be a compact subset of $X$. Let fix $\epsilon\gg \theta$ and $\lambda\in (0,1)$. Now,
  we will consider an $(\epsilon,\lambda)$-cover $\{N_x(\epsilon,\lambda):x\in A\}$. Since $A$
  is compact, there exist $x_1,x_2,\cdots,x_n\in A$ such that
   $A\subseteq \displaystyle{\bigcup_{i=1}^nN_{x_i}(\epsilon,\lambda)}$. Let $x,y\in A$.
   Then there exists $i\in \{1,\cdots n\}$ such that $x\in N_{x_i}(\epsilon,\lambda)$ and exists
   $j\in \{1,\cdots n\}$ such that $y\in N_{x_j}(\epsilon,\lambda)$. Thus we have
   $\F_{x,x_i}(\epsilon)>1-\lambda$ and $\F_{y,x_j}(\epsilon)>1-\lambda$. Now, let
   $m=\min\{\F_{x_i,x_j}(\epsilon):1\leq i,j\leq n\}$. It is obvious that $m>0$ and we have
   $$\F_{x,y}(\epsilon)\geq \F_{x,x_i}(\epsilon)\ast \F_{x_i,x_j}(\epsilon)\ast \F_{x_j,y}(\epsilon)
  \geq (1-\lambda)\ast m\ast(1-\lambda)>1-\delta,$$
   for some $0<\delta<1$. If we take $\epsilon_1=3\epsilon$, we have $\F_{x,y}(\epsilon)>1-\delta$ for all $x,y\in A$.
   Hence we obtain that $A$ is probabilistic semi-bounded set.
\end{proof}
\section{Convex Structure, Normal structure and Strictly Convex Structure on FM-Spaces}
The idea of convex metric spaces was initially introduced by Takahashi \cite{Takahashi}.
Both normed linear spaces and hyperbolic metric spaces fall under this group of metric spaces.
This section aims at giving definitions of strictly convex and normal structures in probabilistic cone metric space.
\begin{definition}
  Let $(X,d)$ be a cone metric space. We say that a cone metric space possesses a Takahashi convex structure
  if there exists a function $W: X^2\times [0,1]\rightarrow X$ which satisfies
  $$d(z,W(x,y,\mu))\leq \mu d(z,x)+(1-\mu)d(z,y),$$
  for all $x,y,z\in X$ and arbitrary $\mu\in [0,1]$. A cone metric space $(X,d)$ with Takahashi's structure is called convex cone metric
  space.
  \end{definition}
  \indent  In this section, we present an extension of Takahashi's notion to the context of a probabilistic cone metric space.
  \begin{definition}
  Let $(X,\F,\ast)$ be a probabilistic cone metric space with continuous $t$-norm $\ast$. A mapping
  $S:X\times X\times [0,1]\rightarrow X$ is said to be a convex structure on $X$ if for
  every $(x,y)\in X\times X$ holds $S(x,y,0)=y$, $S(x,y,1)=1$ and for all $x,y,z\in X$, $\mu\in [0,1]$
  and $\epsilon\gg\theta$
  \begin{equation}\label{Eq.G1}
  \F_{S(x,y,\mu),z}(2\epsilon)\geq \F_{x,z}\left(\frac{\epsilon}{\mu}\right)\ast\F_{x,z}\left(\frac{\epsilon}{1-\mu}\right).
  \end{equation}
\end{definition}
\indent It is easy to see that every cone metric space $(X,d)$ with a convex structure $S$ can be consider as
a probabilistic cone metric space $(X,\F,\ast_{\min})$ (the associated probabilistic cone metric space) with the same function $S$.
A probabilistic cone metric space $(X,\F,\ast)$ with a convex structure $S$ is called a convex probabilistic cone metric space.\\
\indent \indent The terminology used in the sequel are provided here.
\begin{definition}
  A point $x\in A$ is called diametral if
  $\displaystyle{\inf_{y\in A}\sup_{s<t}\F_{x,y}(s)=\delta_A(t)}$ holds for all $t>0$.
\end{definition}
\begin{definition}
Let $(X,\F,\ast)$ be a probabilistic cone metric space with continuous $t$-norm $\ast$ and a convex structure
$S(x,y,\mu)$. A subset $A\subseteq X$ is said to be a convex if for every $x,y\in A$ and $\mu\in [0,1]$
it follows that $S(x,y,\mu)\in A$.
\end{definition}
\begin{lemma}\label{LemmaA}
Let $(X,\F,\ast)$ be a probabilistic cone metric space with continuous $t$-norm $\ast$ and let
$\{K_{\alpha}\}_{\alpha\in\Lambda}$  be a family of convex subsets of $X$. Then the intersection
$K=\bigcap_{\alpha\in\Lambda} K_{\alpha}$ is a convex set.
\end{lemma}
\begin{proof}
  If $x,y\in K$, then $x,y\in K_{\alpha}$ for every $\alpha\in\Lambda$. It follows that $S(x,y,\mu)\in K_{\alpha}$ for
  every $\alpha\in\Lambda$, i.e., $S(x,y,\mu)\in K$, which means that $K$ is convex.
\end{proof}
\begin{definition}
  A convex probabilistic cone metric space $(X,\F,\ast)$ with a convex structure $S:X\times X\times [0,1]\rightarrow X$ and continuous
  $t$-norm $\ast$ will be called strictly convex if, for arbitrary $x,y\in X$ and $\mu\in (0,1)$ the element $z=S(x,y,\mu)$
  is the unique element which satisfies
  \begin{equation}\label{Eq.G2}
 \F_{x,y}\left(\frac{\epsilon}{\mu}\right)=\F_{x,y}(\epsilon),\quad \F_{x,y}\left(\frac{\epsilon}{1-\mu}\right)=\F_{x,y}(\epsilon),
  \end{equation}
for all $\epsilon\gg\theta$.
\end{definition}
\begin{lemma}\label{LemmaB}
  Let $(X,\F,\min)$ be a probabilistic cone metric space with a convex structure $S(x,y,\mu)$ and continuous
  $t$-norm $\ast$. Suppose that for every $\mu\in (0,1)$, $t>0$ and $x,y,z\in X$ hold
  \begin{equation}\label{Eq.G3}
 \F_{S(x,y,\mu),z}(\epsilon)>\min\{\F_{x,z}(\epsilon),\F_{z,y}(\epsilon)\}.
  \end{equation}
  If there exists $z\in X$ such that
  \begin{equation}\label{Eq.G4}
   \F_{S(x,y,\mu),z}(\epsilon)=\min\{\F_{x,z}(\epsilon),\F_{z,y}(\epsilon)\}
  \end{equation}
  is satisfied, for all $\epsilon\gg \theta$, then $S(x,y,\mu)\in\{x,y\}$.
\end{lemma}
\begin{proof}
  Let us assume that (\ref{Eq.G4}) holds for some $z\in X$ and for all $\epsilon\gg \theta$. Since (\ref{Eq.G3})
  holds, it follows that $\mu=0$ or $\mu=1$ and, consequently we have that $S(x,y,0)=y$ or $S(x,y,1)=x$,which proves the
  lemma.
\end{proof}
\begin{lemma}\label{LemmaC}
   Let $(X,\F,\ast)$ be a probabilistic cone metric space  with a convex structure $S(x,y,\mu)$ and continuous
  $t$-norm $\ast$. Then for arbitrary $x,y\in X$, $x\neq y$ there exists $\mu\in (0,1)$ such that
  $S(x,y,\mu)\notin \{x,y\}$.
\end{lemma}
\begin{proof}
  Suppose that for every $\mu\in(0,1)$, it holds that $S(x,y,\mu)\in \{x,y\}$. From (\ref{Eq.G2}) it follows that
  $M(x,y,t)=1$ for all $t>0$ which means that $x=y$ and so the proof is achieved.
\end{proof}
\begin{definition}
  A probabilistic cone metric space $(X,\F,\ast)$ with continuous $t$-norm $\ast$ possesses a normal structure
  if, for every closed , probabilistic semi-bounded and convex set $Y\subset X$, which consists of at least
  two different points, there exists a point $x\in Y$ which is non-diametral, i.e., there exists $t_0>0$
  such that
  $$\delta_Y(t_0)<\displaystyle{\inf_{y\in Y}\sup_{s<t_0}\F_{x,y}(s)}$$
  holds.
\end{definition}
It is obvious that compact and convex sets in convex cone metric space possess a normal structure.
\begin{definition}
  Let $(X,\F,\ast)$ be a convex probabilistic cone metric space with continuous $t$-norm $\ast$ and $Y\subseteq X$.
  The closed convex shell of a set $Y$ denoted by $cov(Y)$, is the intersection of all closed, convex sets that contain $Y$.
\end{definition}
Since the collection of closed, convex sets that contain $Y$ is not empty and $X$ is a member of this collection, it is obvious that the set $cov(Y)$ exists. It follows that this intersection is a convex set from Lemma \ref{LemmaA}.
Additionally, as an intersection of closed sets, this intersection is closed.
\begin{definition}
  Let $(X,\F,\ast)$ be a probabilistic cone metric space with continuous $t$-norm $\ast$ and let
  $f$ be a self-mapping on $X$. We say that $f$ is a non-expansive mapping if
  \begin{equation}\label{EqG5}
     \F_{fx,fy}(\epsilon)\geq \F_{x,y}(\epsilon)
  \end{equation}
  holds for all $x,y\in X$ and $\epsilon\gg \theta$.
\end{definition}
\begin{example} \label{Example-A} Let $X=\mathbb{R}$, $E=\mathbb{R}$ and $P=[0,\infty)$ be a cone.
	Let $S:X\times X\times[0,1]\rightarrow X$ defined by
	$$S(x,y,\mu)=\mu x+(1-\mu)y$$
	for all $x,y\in \mathbb{R}$ and $\mu\in (0,1)$ is a convex structure on probabilistic cone metric space $(\R,\F,\ast)$
	induced by a metric $d(x,y)=|x-y|$ on $X$, where $a\ast b=\min\{a,b\}$ is continuous $t$-norm for $a,b\in [0,1]$
	and
	$$\F_{x,y}(t)=\left\{
	\begin{array}{ll}
	0, & \hbox{if $t\leq d(x,y)$;} \\
	1, & \hbox{if $d(x,y)< t$.}
	\end{array}
	\right.
	$$
	for all $x,y\in X$ and $t\gg \theta$. Let us prove this assertion. Firstly, we have that
	$S(x,y,0)=y$ and $S(x,y,1)=x$ for all $x,y\in X$. Now, let us prove that inequality (\ref{Eq.G1}) is satisfied.
	If we assume that
	$$\min\left\{\F_{x,z}\left(\frac{t}{\mu}\right),\F_{y,z}\left(\frac{t}{1-\mu}\right)\right\}=0$$
	then inequality (\ref{Eq.G1}) is a trivially satisfied because we get $\F_{S(x,y,\mu),z}(2t)\geq 0$. Now we will assume
	that $\F_{x,z}\left(\frac{t}{\mu}\right)=1$ and $\F_{y,z}\left(\frac{t}{1-\mu}\right)=1$. Then we have
	that $\frac{t}{\mu}>d(x,z)$ and $\frac{t}{1-\mu}>d(y,z)$, i.e. $t>\mu d(x,z)$ and $t>(1-\mu)d(y,z)$. Hence, we obtain
	\begin{align*}
	2t &>\mu d(x,z)+(1-\mu)d(y,z)=\mu|x-z|+(1-\mu)|y-z| \\
	&\geq |\mu x-\mu z+(1-\mu)y-(1-\mu)z|\\
	&=|\mu x+(1-\mu)y-z|=d(\mu x+(1-\mu)y,z)\\
	&=d(S(x,y,\mu),z),
	\end{align*}
	that is, $2t-d(S(x,y,\mu),z)>0$ and so $M(S(x,y,\mu),z,2t)=1$. Therefore, inequality (\ref{Eq.G1}) holds
	for all $x,y,z\in X$ and $t\gg \theta$.
\end{example}
\indent It is easy to see that every cone metric space $(X,d)$ with a convex structure $S$ can be consider as
a probabilistic cone metric space $(X,\F,\ast)$ (the associated probabilistic metric space) with the same function $S$.
A probabilistic cone metric space $(X,\F,\ast)$ with a convex structure $S$ is called a convex probabilistic cone metric space.\\

\section{Common Fixed Point Theorems}
The goal of this section is to show that two self-mappings formed on strictly convex probabilistic cone metric space have a common fixed point theorem. Topological approaches for characterisation spaces with nondeterministic distance will be utilized to prove the primary conclusion.
\begin{lemma}\label{Lemma4.1} Let $K\subseteq X$  be a non-empty, convex, and compact subset of $X$. Also, let $(X,F,\ast)$ be a strictly convex probabilistic cone metric space with continuous $t$-norm $\ast$ and convex structure $S(x,y,\mu)$ fulfilling (\ref{Eq.G3}). Then, $K$ has a normal structure.
\end{lemma}
\begin{proof}
 Let's say that $K$ doesn't have a normal structure. Then, for any $x$ in $Y$, there exists a closed, probabilistic, semi-bounded, convex subset $Y\subset K$ that contains at least two distinct points such that $Y$ does not contain a non-diametral point, i.e.,$$\displaystyle{\inf_{y\in Y}\sup_{s<t}\F_{x,y}(s)=\delta_{Y}(t)}$$
  for every $x\in Y$.
Since condition (\ref{Eq.G3}) is satisfied and $X$ is strictly convex, the claims made by lemma \ref{LemmaB} and Lemma \ref{LemmaC} are true. Assume that $x_1$ and $x_2$ be arbitrary points in $Y$. The fact that there is a $\mu_0\in (0,1)$ such that $S(x_1,x_2,\mu_0)\notin \{x_1,x_2\}$ exists is inferred from Lemma \ref{LemmaC}. The fact that $Y$ is a convex set means that $S(x_1,x_2,\mu_0)\in Y$.
Given that $Y$ is a closed subset of $K$, which is a compact set, $Y$ is also compact. Since $\delta_Y(t)=\inf_{y\in Y}\sup_{s<t}\F_{y,S(x_1,x_2,\mu)}(s)$ is a left-continuous function on the compact set $Y$ for every $t>0$, there exists $x_3,x_4\in Y$ such that $\sup_{s<t}\F_{x_3,S(x_1,x_2,\mu_0)}(s)=\delta_Y(t)$. Given that  $\F_{x,y}(.)$  is a non-decreasing left-continuous function and Lemma \ref{LemmaB} it follows that
  \begin{align}\label{EqG6}
    \delta_Y(t) &=\sup_{s<t}\F_{x_3,S(x_1,x_2,\mu_0)}(s)=\F_{x_3,S(x_1,x_2,\mu_0)}(t)\nonumber \\
     &> \min\{\F_{x_3,x_1}(t),\F_{x_3,x_2}(t)\}\\
     &=\min\{\sup_{s<t}\F_{x_3,x_1}(s),\sup_{s<t}\F_{x_3,x_2}(s)\}\geq \delta_Y(t).\nonumber
  \end{align}
  It follows from the previous inequality that $\delta_Y(t)>\delta_Y(t)$, which is a contradiction. The evidence is now complete.
\end{proof}
\begin{lemma}\label{Lemma4.2}
  Let $(X,\F,\ast)$ be a convex probabilistic cone metric space with a convex structure $S(x,y,\mu)$ satisfying the
  condition (\ref{Eq.G3}). Then the closed $(\epsilon,\lambda)$-neigborhoods $N_{x}[\epsilon,\lambda]$ are
  convex sets.
\end{lemma}
\begin{proof}
  Let $a,b\in N_{x}[\epsilon,\lambda]$ be arbitrary points. This implies that
  $\F_{a,x}(\epsilon)\geq 1-\lambda$ and $\F_{b,x}(\epsilon)\geq 1-\lambda$ for all
  $\epsilon>0$. We shall prove that $\F_{S(a,b,\mu),x}(\epsilon)\geq 1-\lambda$ for
  all $\epsilon>0$, i.e., $S(a,b,\mu)\in N_{x}[\epsilon,\lambda]$. Indeed, for
  $\mu\in (0,1)$, from (\ref{Eq.G3}), we have that
  $$\F_{S(a,b,\mu),x}(\epsilon)>\min\{\F_{a,x}(\epsilon),\F_{b,x}(\epsilon)\}\geq \min\{1-\lambda,1-\lambda\}=1-\lambda.$$
  For $\mu=0$ or $\mu=1$. It follows that $S(a,b,0)=b$ or $S(a,b,1)=a$ belongs to $N_{x}[\epsilon,\lambda]$.
\end{proof}
\begin{example} \label{Example-B}
	Let $E=\R$. Then $P=\{x\in\R:x\geq 0\}$ is a normal cone with normal constant $K=1$.
	Let $X=\R$, $a\ast b=\min\{a,b\}$ and $\F:X^2\times int(P)\rightarrow [0,1]$ defined by
	$$\F_{x,y}(t)=H(t-d(x,y))$$
	for all $x,y\in X$ and $t\gg \theta$, where $d:X^2\rightarrow E$ is a cone metric space defined by $d(x,y)=|x-y|$ and
	$$H(t)=\left\{
	\begin{array}{ll}
	0, & \hbox{if $d(x,y)\leq t$;} \\
	1, & \hbox{if $d(x,y)>t$.}
	\end{array}
	\right.
	$$
	Then $(X,\F,\ast)$ is a probabilistic  cone metric spaces.\\
	The mapping $S:\R\times\R\times [0,1]\rightarrow [0,1]$ defined by
	$$S(x,y,\mu)=\mu x+(1-\mu)y$$
	for all $x,y\in\R$ and $\mu\in (0,1)$ is a convex structure on probabilistic cone metric space.
	For arbitrary $x,y\in X$ and $\mu\in (0,1)$ the element $z=S(x,y,\mu)=\mu x+(1-\mu)y$ is the unique
	element which satisfies
	\begin{align*}
	\F_{z,x} (t)&=H(t-d(x,z))=H(t-d(x-\mu x,z-\mu x)) \\
	&=H\left(\frac{t}{1-\mu}-\frac{d((1-\mu)x,z-\mu x)}{1-\mu}\right)\\
	&=H\left(\frac{t}{1-\mu}-d\left(x,\frac{z}{1-\mu}-\frac{\mu x}{1-\mu}\right)\right)\\
	&=H\left(\frac{t}{1-\mu}-d(x,y)\right)\\
	&=\F_{x,y}\left(\frac{t}{1-\mu}\right).
	\end{align*}
	In a similar way it can proved that the second equality in (\ref{Eq.G2})
	is satisfied. Hence we obtained that the probabilistic cone metric space is strictly convex
	with a given convex structure $S(x,y,\mu)$.\\
	\indent On the other hand, we have that
	$$d(\mu x+(1-\mu)y,z)<\max\{d(x,z),d(y,z)\}$$
	is satisfied for all $\mu\in (0,1)$, and it follows that
	\begin{align*}
	\F_{S(x,y,\mu),z} (t) &=H(t-d(S(x,y,\mu),z)) \\
	&>H(t-\max\{d(x,z),d(y,z)\})\\
	&=\min\{H(t-d(x,z)),H(t-d(y,z))\}\\
	&=\min\{\F_{x,z}(T),\F_{y,z}(t)\}
	\end{align*}
	holds, that is, condition (\ref{Eq.G3}) is satisfied.
\end{example}
Now we give the main theorem in this paper.
\begin{theorem}\label{Theorem4.3}
Let $(X,\F,\ast)$ be a convex probabilistic cone metric space with a convex structure $S(x,y,\mu)$ satisfying the
  condition (\ref{Eq.G3}) and let $E\subseteq X$ be a non-empty, convex and compact subset of $X$. Let $f$ and $g$ be a
  self-mappings on $E$, $g(E)\cap f(E)\subseteq E$, satisfying the conditions
  \begin{equation}\label{EqG7}
    \F_{f(x),g(y)}(t)\geq \F_{x,y}(t)
  \end{equation}
  for all $x,y\in E$, $x\neq y$ and for every $t\gg \theta$. Then $f$ and $g$ have at least one common fixed point on $E$.
\end{theorem}
\begin{proof}
  First, note that the collection $\mathfrak{C}$ of all non-empty, closed, convex sets $E_{\alpha}\subseteq E$
  such that $g(E_{\alpha})\cap f(E_{\alpha})\subseteq E_{\alpha}$ is non-empty, because $E\subseteq \mathfrak{C}$.
  Indeed, $E$ is closed set because the fact that it is compact set in a Hausdorff pace and it is satisfied that
  $g(E)\cap f(E)\subseteq E$. If we order this collection with inclusion, then $(\mathfrak{C},\subseteq)$ is a partially ordered set.
  Let $E_{\alpha}:\alpha\in\Lambda$ be an arbitrary chain of this family. Then
  the set $\displaystyle{\bigcap_{\alpha\in \Lambda}E_{\alpha}}$ is nonempty, closed, convex subset of $E$, which is a lower bound of this chain.
  Indeed, let us assume that $\displaystyle{\bigcap_{\alpha\in \Lambda}E_{\alpha}=\emptyset}$. Then, from Lemma \ref{Lemma2.9}, it follows that there exists a finite sub-collection $E_{\alpha_1}\supseteq E_{\alpha_2}\supseteq \cdots \supseteq E_{\alpha_n}$ of
  the chain $\{E_{\alpha}:\alpha\in\Lambda$ which has an empty intersection, which is impossible, since this intersection is
  $E_{\alpha_n}\neq \emptyset$.  By Zorn's Lemma, it follows that there exists a minimal element $E_0$ of the collection $\mathfrak{C}$ such
  that $g(E_0)\cap f(E_0)\subseteq E_0$. We will prove that $E_0$ consists of only one point and since
  $g\cap f:E_0\rightarrow E_0$ this will mean that $g$ and $f$ have a common fixed point.\\
  \indent Let us assume that $E_0$ contains at least two different points.It follows from Lemma \ref{Lemma4.1} that $E$ possesses
  a normal structure. It follows From Theorem \ref{Theorem2.13} that $E_0$ is probabilistic semi-bounded set. Since
  $E_0$ is closed and convex, it follows that there exists some non-diametral point $x_0\in E_0$, i.e.,
  there exists $t_0\gg \theta$ such that the following inequality holds:
  \begin{equation}\label{EqG8}
    \displaystyle{\inf_{y\in E_0}\sup_{s<t}\F_{x_0,y}(s)>\delta_{E_0}(t_0)}.
  \end{equation}
  Denote $1-\rho:=\displaystyle{\inf_{y\in E_0}\sup_{s<t}\F_{x_0,y}(s)}.$ Let us denote with $E_1$
  the closed convex shell of the set $g(E_0)\cap f(E_0)$. Since $g(E_0)\cap f(E_0)\subseteq E_0$, it holds that
  $E_1=cov(g(E_0)\cap f(E_0))=\overline{cov(g(E_0)\cap f(E_0))}\subseteq \overline{cov(E_0)}=\overline{E_0}=E_0$.
   Therefore, $E_1\subseteq E_0$ and it follows that $g(E_1)\cap f(E_1)\subseteq g(E_0)\cap f(E_0)\subseteq \overline{cov(g(E_0)\cap f(E_0))}=E_1,$
   i.e., $g(E_1)\cap f(E_1)\subseteq E_1.$ This means that $E_1\in \mathfrak{C}$, and since $E_0$ is the minimal element, we have that
   $E_1=E_0$.\\
   \indent If the inequality (\ref{EqG8}) holds, i.e., if $1-\rho>\delta_{E_0}(t_0)$, let us define sets
  $$A:=\displaystyle{\left(\bigcap_{y\in E_0}N_{y}[\rho,t_0]\right)\bigcap E_0}\,\,\text{and}\,\,
            A_1:=\displaystyle{\left(\bigcap_{y\in g(E_0)\cap f(E_0)}N_{y}[\rho,t_0]\right)\bigcap E_0}.$$
\indent The set $A$ is nonempty since $x_0\in A$. Indeed , from inequality (\ref{EqG8}), it follows that
$\F_{x_0,y}(t_0)\geq 1-\rho$.  From the previous we conclude that $x_0\in N_{y}[\rho,t_0]$ for all $y\in E_0$. Consequently,
$x_0$ belongs to $A$. We will show that $A=A_1$. Since $g(E_0)\cap f(E_0)\subseteq E_0$, we have $A\subset A_1$.\\
\indent Now, let $z\in A_1$, we will prove that $z\in A$. Since $z\in A_1$, for arbitrary $y\in g(E_0)\cap f(E_0)$,
it holds that $\F_{y,z}(t_0)\geq 1-\rho$, i.e., $y\in N_{y}[\rho,t_0]$. Since $y$ is an arbitrary point from
$g(E_0)\cap f(E_0)$, it follows that $g(E_0)\cap f(E_0)\subseteq N_{y}[\rho,t_0]$. Because of the fact that $N_{z}[\rho,t_0]$
 is a closed and convex set which contains $g(E_0)\cap f(E_0)$, we conclude that
 $$E_1\subseteq \overline{cov(g(E_0)\cap f(E_0))}\subseteq N_{z}[\rho,t_0] $$
 holds. Since $E_0=E_1$, it follows that $E_0\subseteq N_{z}[\rho,t_0]$. From last we have that for every $y\in E_0$, it holds that
 $z\in N_{y}[\rho,t_0]$, which means that $A_1\subseteq A$, i.e., $A=A_1$.\\
 \indent We will show that $A\in \mathfrak{C}$. Note that $A$ is closed as an intersection of closed sets. From Lemma \ref{LemmaA}
 and Lemma \ref{Lemma4.2} that $A$ is a convex set. Let us prove that $g(A)\cap f(A)\subseteq A$. Let $z\in A$ and $y\in g(E_0)\cap f(E_0)$.
 Then there exists $x\in E_0$ such that $y=f(x)$ and $y=f(x)$. Applying inequality (\ref{EqG7}) for $t=t_0$, we have
 $$\F_{f(z),y}(t_0)=\F_{f(z),g(x)}(t_0)\geq\F_{z,x}(t_0)\geq 1-\rho.$$
 This means that $f(z)\in A$. Since $z$ is an arbitrary point from $A$, we obtain $f(A)\subseteq A$.\\
 \indent On the other hand, we have
 $$\F_{g(z),y}(t_0)=\F_{g(z),f(x)}(t_0)\geq \F_{z,x}(t_0)\geq 1-\rho.$$
 This means that $g(z)\in A$. Since $z$ is an arbitrary point from $A$, we obtain $g(A)\subseteq A$.
 Finally, we obtain $g(A)\cap f(A)\subseteq A$.\\
 \indent Since $A\subseteq E_0$ and $E_0$ is the minimal element of the collection $\mathfrak{C}$, it follows that $A=E_0$. Now, we
 have that $\delta_A(t_0)\geq 1-\rho\geq \delta_{E_0}(t_0)$. This is a contradiction with $A=E_0$, i.e.,
 the assumption that $E_0$ contains at least two different points is wrong, which means that $E_0$ contains only one point which is a common fixed point of the mappings $g$ and $f$. This achieves the proof.
\end{proof}
\begin{example}\label{Example-C}
	Let  $X=\R$. Then $P=\{x\in\R: x\geq 0\}$ is a normal cone with normal constant $K=1$. Let $E=[0,1]$, $a\ast b=\min\{a,b\}$
	for all $a,b\in X$. Define $\F:E^2\times int(P)\rightarrow [0,1]$ by
	$$\F_{x,y}(t)=\frac{t}{t+d(x,y)}$$
	for all $x,y\in E$, $d(x,y)=|x-y|$ and $t\gg \theta$. Then $(X,\F,\ast)$ is a probabilistic cone metric space.\\
	Let $f,g:E\rightarrow E$ be a mapping defined by $fx=gx=\frac{x}{2}$
	and  we show that $f$ and $g$ satisfies (\ref{EqG7}). We have
	\begin{align*}
	\F_{fx,gy}(t) &=\frac{t}{t+|fx-gy|}  \\
	&=\frac{t}{t+|\frac{x}{2}-\frac{y}{2}|} \\
	&\geq \frac{t}{t+|\frac{x}{2}-\frac{y}{2}|}\\
	&\geq \frac{t}{t+\frac{1}{2}d(x,y)}\\
	&\geq \frac{t}{t+d(x,y)}\\.
    &\geq \F_{x,y}(t)
	\end{align*}
	This shows that (\ref{EqG7}) and all conditions of Theorem \ref{Theorem4.3} hold and $f$ and $g$ have a common point in $X$. The common point is 0.
\end{example}
\begin{corollary}\label{lastcor}
  Let $(X,\F,\ast)$ be a strictly convex probabilistic cone metric space with a convex structure $S(x,y,\mu)$ satisfying the
  condition (\ref{Eq.G3}) and let $E\subseteq X$ be a non-empty, convex and compact subset of $X$. Let $f$  be a non-expansive
  self-mappings on $E$. Then $f$ has at least one  fixed point on $E$.
\end{corollary}
\begin{proof}
  Putting in the Theorem \ref{Theorem4.3} that $f=g$, we have that the mapping $f$ is a self-mapping on $K$ and in this case,
  from conditions (\ref{EqG7}) and (\ref{EqG8}) we obtain that mapping $f$ is non-expansive on $K$, and hence the conclusion
  holds.
\end{proof}
\begin{example}  Let $X=\mathbb{R}, E=[0,1]$ and $P=[0,\infty)$ and let $d$ defined by
	$d(x,y)=|x-y|$ and let $f:E\rightarrow E$ defined by $f(x)=\frac{x^2}{3}+\frac{1}{2}$.\\
	\indent Let us show that all conditions of Corollary \ref{lastcor} are satisfied. From Example
	\ref{Example-A} and Example \ref{Example-B} we have that $(X,\F,\min)$ is a strictly convex probabilistic
	cone metric space satisfies condition (\ref{Eq.G3}) with a convex structure
	$S(x,y,\mu)=\mu x+(1-\mu)y$ for all $x,y\in X$ and $\mu\in (0,1)$. Notice that $\frac{1}{3}|x^2-y^2|\leq |x-y|$ for all $x,y\in E$,
because $\frac{1}{3}|x+y|\leq 1$ for all $x,y\in E$.
 The mapping  $f$
	is a non-expansive self-mapping on $E$. Indeed, if $x,y\in E$, then
		\begin{align*}
	\F_{fx,fy}(t) &=H(t-d(fx,fy))=H(t-|fx-fy|) \\
	&=H\left(t-\left|\left(\frac{x^2}{3}+\frac{1}{2}\right)-\left(\frac{y^2}{3}+\frac{1}{2}\right)\right|\right)\\
          & = H\left(t-\frac{1}{3}|x^2-y^2|\right)\\
	&\geq H(t-|x-y|)=H(t-d(x,y))=\F_{x,y}(t),
	\end{align*}
It is obvious that $E$ is nonempty and compact set.
	that is, condition (\ref{EqG7}) is satisfied for every $x,y\in E$. Since all conditions of Corollary \ref{lastcor}
	are satisfied and so we obtain that $f$ has a unique fixed point $x=\frac{3-\sqrt{3}}{2}\in E$.
\end{example}
{\bf Author Contributions:}   The author have read and agreed to the published version of the manuscript.\\
{\bf Funding:} No funding is applicable\\
{\bf Institutional Review Board Statement:} Not applicable.\\
{\bf Informed Consent Statement:} Not applicable.
{\bf Data Availability Statement:} Not applicable.\\
{\bf Conflicts of Interest:} The authors declare no conflict of interest.

\bibliographystyle{unsrtnat}
\bibliography{references}  






\end{document}